\documentstyle[11pt]{article}

\newcommand{\bd}{
\begin{document}}  
\newcommand{\ed}{\end{document}}
\newcommand{\bc}{\begin{center}}
\newcommand{\ec}{\end{center}}
\newcommand{\bq}{\begin{quote}}
\newcommand{\eq}{\end{quote}}
\newcommand{\lb}{\linebreak}
\renewcommand{\pb}{\pagebreak}
\newcommand{\mb}{\makebox}
\newcommand{\lt}{\left}
\newcommand{\rt}{\right}

\newcommand{\vs}{\vspace}
\newcommand{\hf}{\hspace*{\fill}}
\newcommand{\vf}{\vspace*{\fill}}
\newcommand{\beqa}{\begin{eqnarray*}}
\newcommand{\eeqa}{\end{eqnarray*}}
\newcommand{\beqn}{\begin{eqnarray}}
\newcommand{\eeqn}{\end{eqnarray}}
\newcommand{\bbibl}{\begin{thebibliography}{99}}
\newcommand{\ebibl}{\end{thebibliography}}
\newcommand{\bm}[1]{\mb{\boldmath ${#1}$}}
\newcommand{\fb}[1]{\lt[{#1}\rt]}
\newcommand{\ot}{\otimes}
\newcommand{\nn}{\nonumber}

\newcommand{\R}{\mb{$I\!\!R$}}
\newcommand{\N}{I\!\!N}
\newcommand{\C}{{\cal C}}
\newcommand{\M}{{\cal M}}
\newcommand{\A}{{\cal A}}
\renewcommand{\P}{{\cal P}}
\renewcommand{\S}{{\cal S}}

\newcommand{\ba}{\begin{array}}
\newcommand{\ea}{\end{array}}
\newcommand{\et}[2]{#1_{1}, #1_{2}, \ld , #1_{#2}}
\newcommand{\eti}[2]{#1_{i1}, #1_{i2}, \ld , #1_{i#2}}
\newcommand{\alt}[4]{\lt\{ \ba{ll}#1 & \mb{ if \quad}#2 \\ #3 & \mb{ if
               \quad}#4 \ea \rt.} 
\newcommand{\alter}[6]{\lt\{ \ba{ll}#1 & \mb{ if \quad}#2 \\ #3 & \mb{ if 
               \quad} #4 \\ #5 & \mb{ if \quad}#6 \ea \rt.}
\newcommand{\alto}[3]{\lt\{ \ba {ll}#1 & \mb{ if \quad}#2 \\ #3 &
               \mb{ otherwise} \ea \rt.} 
\newcommand{\altero}[5]{\lt\{ \ba {ll}#1 & \mb{ if \quad}#2 \\ #3 &
               \mb{ if \quad} #4 \\ #5 & \mb{ otherwise} \ea \rt.} 

\newcommand{\es}{\emptyset}
\newcommand{\ci}{\subseteq}
\newcommand{\cs}{\supseteq}

\renewcommand{\u}{\cup}
\renewcommand{\i}{\cap}
\newcommand{\bu}{\bigcup}
\newcommand{\bi}{\bigcap}

\newcommand{\ra}{\rightarrow}
\newcommand{\Ra}{\Rightarrow}
\newcommand{\Lra}{\Leftrightarrow}
\newcommand{\lgra}{\longrightarrow}
\newcommand{\Lgra}{\Longrightarrow}
\newcommand{\lglra}{\longleftrightarrow}
\newcommand{\Lglra}{\Longleftrightarrow}

\renewcommand{\a}{\alpha}
\renewcommand{\b}{\beta}
\newcommand{\G}{\Gamma}
\newcommand{\g}{\gamma}
\renewcommand{\d}{\delta}
\newcommand{\e}{\varepsilon}
\newcommand{\h}{\eta}
\newcommand{\la}{\lambda}
\newcommand{\La}{\Lambda}
\newcommand{\m}{\mu}
\newcommand{\p}{\pi}
\newcommand{\s}{\sigma}
\newcommand{\Si}{\Sigma}
\newcommand{\ta}{\tau}
\newcommand{\ph}{\phi}
\newcommand{\Ph}{\Phi}
\renewcommand{\c}{\chi}
\newcommand{\om}{\omega}
\newcommand{\Om}{\Omega}

\newcommand{\pf}{\noindent {\bf Proof :} }
\newcommand{\rec}[1]{\frac{1}{#1}}
\newcommand{\f}{\frac}
\newcommand{\sm}[2]{\sum_{#1}^{#2}}
\newcommand{\ld}{\ldots}
\newcommand{\ov}{\overline}
\newcommand{\un}{\underline}
\newcommand{\iy}{\infty}
\newcommand{\qed}{\hf\rule{6pt}{6pt} \vs{.5\smallskipamount}}
\newcommand{\wt}{\widetilde}
\newcommand{\ds}{\displaystyle}

\newcounter{cnt1}
\newcounter{cnt2}
\newcounter{cnt3}
\newcommand{\blr}{\begin{list}{$($\roman{cnt1}$)$} {\usecounter{cnt1}
		\setlength{\topsep}{0pt} \setlength{\itemsep}{0pt}}}
\newcommand{\bla}{\begin{list}{$($\alph{cnt2}$)$} {\usecounter{cnt2}
		\setlength{\topsep}{0pt} \setlength{\itemsep}{0pt}}}
\newcommand{\bln}{\begin{list}{$($\arabic{cnt3}$)$} {\usecounter{cnt3}
                \setlength{\topsep}{0pt} \setlength{\itemsep}{0pt}}}
\newcommand{\el}{\end{list}}

\newtheorem{thm}{Theorem}
\newtheorem{lem}[thm]{Lemma}
\newtheorem{cor}[thm]{Corollary}
\newtheorem{ex}{Example}
\newtheorem{Q}{Question}
\newtheorem{F}[thm]{Fact}
\newtheorem{Def}{Definition}
\newtheorem{prop}[thm]{Proposition}
\newtheorem{rem}{Remark}
\newcommand{\Rem}{\begin{rem} \rm}
\newcommand{\bdfn}{\begin{Def} \rm}
\newcommand{\edfn}{\end{Def}}
\newcommand{\bft}{\begin{F} \rm}
\newcommand{\eft}{\end{F}}
\newcommand{\TFAE}{the following are equivalent : }

\renewcommand{\baselinestretch}{1.5}

\title{\bf On Nicely Smooth Banach Spaces}
\author{\bf Pradipta Bandyopadhyay \& Sudeshna Basu\\
\sc Stat--Math Division, Indian Statistical Institute\\ \sc 203,
B.\ T.\ Road, Calcutta 700 035, India\\ e-mail~:
pradipta@isical.ernet.in \& res9415@isical.ernet.in}
\date{}
\sloppy
\bd

\maketitle

\thispagestyle{empty}
\vfill

\begin{abstract}
In this work, we obtain some necessary and some sufficient
conditions for a space to be nicely smooth, and show that they
are equivalent for separable or Asplund spaces. We obtain a
sufficient condition for the Ball Generated Property (BGP), and
conclude that Property $(II)$ implies the BGP, which, in turn,
implies the space is nicely smooth. We show that the class of
nicely smooth spaces is stable under $c_o$ and $\ell_p$ sums and
also under finite $\ell_1$ sums; that being nicely smooth is not
a three space property; and that the Bochner $L_p$ spaces are
nicely smooth if and only if $X$ is both nicely smooth and
Asplund. A striking result obtained is that every equivalent
renorming of a space is nicely smooth if and only if it is
reflexive.
\end{abstract}
\vfill

\noindent \hrulefill\
\begin{description} \small
\item[AMS Subject Classification (1990) :] 46B20, 46B22.
\item[Keywords and Phrases :] Nicely smooth spaces, Ball
Generated Property, Property $(II)$. 
\end{description}
\pb

We work with {\em real}\/ Banach spaces. We will denote by
$B(X)$, $S(X)$ and $B[x, r]$ respectively the closed unit ball,
the unit sphere and the closed ball of radius $r>0$ around $x
\in X$.

\bdfn
We say $A \ci B(X^*)$ is a norming set for $X$ if $\|x\| =
\sup\{x^*(x) : x^* \in A\}$, for all $x \in X$. A subspace $F
\ci X^*$ is a norming subspace if $B(F)$ is a norming set for
$X$.

A Banach space $X$ is {\em nicely smooth}\/ if $X^*$ contains no
proper norming subspace; has the {\em Ball Generated Property}\/
(BGP) if every closed bounded convex set in $X$ is
ball-generated, i.e., intersection of finite union of balls; has
{\em Property $(II)$}\/ if every closed bounded convex set in
$X$ is the intersection of closed convex hull of finite union of
balls.
\edfn

In this work, we obtain a sufficient condition for the BGP
weaker than those considered earlier. We conclude that Property
$(II)$ implies the BGP, which, in turn, implies the space is
nicely smooth. We obtain a few more conditions---some necessary,
some sufficient---for a space to be nicely smooth, and show that
they are equivalent for separable or Asplund spaces. We show
that the class of nicely smooth spaces is stable under $c_o$ and
$\ell_p$ sums $(1<p<\iy)$ and also under finite $\ell_1$ sums.
We show that the Bochner $L_p$ spaces $(1<p<\iy)$ are nicely
smooth if and only if $X$ is both nicely smooth and Asplund. A
striking result obtained is that every equivalent renorming of a
space is nicely smooth if and only if it is reflexive. This
significantly strengthens all known results of similar type and
has a surprisingly elementary proof. Some more stability results
are also obtained.

Notice that if a separable space is nicely smooth then it has
separable dual. And a dual space is nicely smooth if and only if
it is reflexive.

For a Banach space $X$, let us denote by $C_X$, the set
$\{x^{**} \in X^{**} : \|x^{**} + \hat{x}\| \geq \|x\|$ for all
$x \in X \}$. Then by \cite[Lemma I.1]{G4}, $x^{**} \in C_X$ if
and only if $ker x^{**}$ is a norming subspace of $X^*$.

\begin{prop} For a Banach space $X$, \TFAE
\bla
\item $X$ is nicely smooth.
\item $C_X = \{0\}$.
\item For all $x^{**} \in X^{**}$,
\[\bi_{x \in X} B[x, \|x^{**} - x\|] = \{x^{**}\}\]
\item Every norming set $A \ci B(X^*)$ separates points of
$X^{**}$.
\el
\end{prop}

\pf Equivalence of (a) and (c) is \cite[Lemma 2.4]{GS}. The rest
are easy. \qed

We now identify some necessary and some sufficient conditions
for a space to be nicely smooth.

For $x\in S(X)$, let $D(x) =\{f\in S(X^*) : f(x) = 1\}$. The set
valued map $D$ is called the duality map and any selection of
$D$ is called a support mapping.

\begin{prop} \label{p1}
For a Banach space $X$, consider the following statements~:
\bla
\item $X^*$ is the closed linear span of the w*-weak PCs of
$B(X^*)$. 
\item Any two distinct points in $X^{**}$ are separated by
disjoint closed balls having centre in $X$.
\item $X$ is nicely smooth.
\item \label{d} For every norm dense set $A \ci S(X)$ and every
support mapping $\phi$, the set $\phi(A)$ separates points of
$X^{**}$.
\el
Then $(a) \Ra (b) \Ra (c) \Ra (d)$.
\end{prop}

\pf $(a) \Ra (b)$ follows from \cite[Theorem 3.1]{CL}. The rest
are easy. \qed

\Rem \label{r1} If we assume in addition that the w*-weak PCs of
$B(X^*)$ form a norming set, then $(a)$, $(b)$ and $(c)$ are
equivalent. And under the stronger assumption that the set
\[\{x \in S(X) : D(x) \mb{ intersects the w*-weak PCs of }
B(X^*)\}\] 
is dense in $S(X)$, all the statements are equivalent.

Can any of the implications be reversed in general?
\end{rem}

The following characterization of the BGP follows from the fact
that $X$ has the BGP if and only if every $x^* \in X^*$ is
ball-continuous on $B(X)$ \cite[Theorem~8.3]{GK} and
\cite[Theorem~1]{CL1}.

\begin{thm} \label{BGP}
X has the BGP if and only if for every $x^* \in X^*$ and $\e >
0$, there exists w*-slices $\et{S}{n}$ of $B(X^*)$ such that for
any $(\et{x^*}{n}) \in \prod_{i=1}^n S_i$, there are {\em
nonnegative} scalars $\et{a}{n}$ such that $\|x^* - \sm{i=1}{n}
a_i x^*_i\| \leq \e$.
\end{thm}

\bdfn
A point $x^*_o$ in a convex set $K \ci X^*$ is called a w*-SCS
point of $K$, if for every $\e > 0$, there exist w*-slices
$\et{S}{n}$ of $K$, and a {\em convex} combination $S =
\sm{i=1}{n} \la_i S_i$ such that $x^*_o \in S$ and diam$(S) <
\e$.
\edfn

\begin{prop} \label{scs}
If $X^*$ is the closed linear span of the w*-SCS points of
$B(X^*)$, then $X$ has the BGP.
\end{prop}

\pf Let $x^* \in X^*$ and $\e > 0$. Since the set of w*-SCS
points of $B(X^*)$ is symmetric and spans $X^*$, there exist
w*-SCS points $\et{x^*}{n}$ of $B(X^*)$, and positive scalars
$\et{a}{n}$ such that $\|x^* - \sm{i=1}{n} a_i x^*_i\| \leq
\e/2$. By definition of w*-SCS points, for each $i =1, 2, \ld,
n$, there exist w*-slices $\eti{S}{n_i}$ of $B(X^*)$, and a
convex combination $S_i = \sm{k=1}{n_i} \la_{ik} S_{ik}$ such
that $x^*_i \in S_i$ and diam$(S_i) < \e/(2 \sm{i=1}{n} a_i)$.
Now, for any $(x^*_{ik}) \in \prod_{i=1}^n \prod_{k=1}^{n_i}
S_{ik}$,
\beqa 
\|x^* - \sm{i=1}{n} \sm{k=1}{n_i} a_i \la_{ik} x^*_{ik}\| & \leq
& \|x^* - \sm{i=1}{n} a_i x^*_i\| + \sm{i=1}{n} a_i \|x^*_i -
\sm{k=1}{n} \la_{ik} x^*_{ik}\| \\ & \leq & \e/2 + \sm{i=1}{n}
a_i {\rm diam}(S_i) \leq \e.
\eeqa

Hence by Theorem~\ref{BGP}, $X$ has the BGP. \qed

\Rem This gives a weaker sufficient condition for the BGP than
ones discussed in \cite[Theorem~7]{CHL}. See Corollary~\ref{ns}
below. 
\end{rem}

\begin{cor} 
Property $(II)$ implies the BGP, which, in turn, implies nicely
smooth.
\end{cor}

\pf Recall that $X$ has Property $(II)$ if and only if w*-PCs of
$B(X^*)$ are norm dense in $S(X^*)$ \cite{CL}, and that a w*-PC
is necessarily a w*-SCS point (this follows from Bourgain's
Lemma, see e.g., \cite[Lemma 1.5]{R}). Thus, Property $(II)$
implies the BGP.

That the BGP implies nicely smooth is proved in \cite{GK}. But
here is an elementary proof.

Let $F$ be a norming subspace of $X^*$. Then $B(X)$ is $\s(X,
F)$-closed, so that every ball-generated set is also $\s(X,
F)$-closed. But if every closed bounded convex set is $\s(X,
F)$-closed, then $F = X^*$. \qed

\begin{cor} \label{ns} If $X$ is an Asplund space (or,
separable), \TFAE
\bla
\item $X^*$ is the closed linear span of the w*-strongly exposed
points of $B(X^*)$.
\item $X^*$ is the closed linear span of the w*-denting points
of $B(X^*)$.
\item $X^*$ is the closed linear span of the w*-SCS points of
$B(X^*)$.
\item $X$ has the BGP.
\item $X$ is nicely smooth, and all the conditions of
Proposition~\ref{p1} are equivalent.
\el
\end{cor}

\pf Clearly, $(a) \Ra (b) \Ra (c) \Ra (d) \Ra$ nicely smooth,
and $(b) \Ra$ Proposition~\ref{p1} (a). Now if $X$ is Asplund
(if $X$ is separable, Proposition~\ref{p1} (d) implies $X^*$ is
separable), then for $A = \{x \in S(X) :$ the norm is Fr\'echet
differentiable at $x\}$, and any support mapping $\phi$,
$\phi(A) = \{$w*-strongly exposed points of $B(X^*)\}$. Hence,
Proposition~\ref{p1} $(d) \Ra (a)$. \qed

\Rem In this case, the assumptions of Remark~\ref{r1} are
satisfied.
\end{rem}

The following result is immediate from \cite[Theorem~12]{CHL} or
\cite[Theorem~2.4 and 2.5]{GK}.

\begin{prop} \label{three}
If every separable subspace of $X$ is nicely smooth, then $X$
has the BGP, and hence, is nicely smooth.
\end {prop}

\bdfn
A Banach space $X$ is said to satisfy finite intersection
property (FIP) if every family of closed balls in $X$ with empty
intersection contains a finite subfamily with empty
intersection.
\edfn

It is well known that any dual space and its 1-complemented
subspaces have FIP.

\begin{thm} \label{four}
$X$ is nicely smooth with FIP if and only if $X$ is reflexive.
\end{thm}

\pf
Sufficiency is obvious.

For necessity, recall from \cite[Theorem~2.8]{GK} that $X$ has
FIP if and only if $X^{**} = X + C_X$. Since $X$ is nicely
smooth, $C_X = \{0\}$ and consequently $X$ is reflexive.
\qed

\Rem Since Hahn-Banach smooth spaces (spaces with Property
$(II)$) are nicely smooth, Theorem~2.11 (Theorem~3.4) of
\cite{br} follow as immediate corollaries with much simpler
proof.
\end{rem}

\begin{thm}
A Banach space $X$ is reflexive if and only if every equivalent
renorming is nicely smooth.
\end{thm}

\pf The converse being trivial, suppose $X$ is not reflexive.
Let $x^{**} \in X^{**} \setminus X$ and let $F = \{x^* \in X^* :
x^{**}(x^*) = 0\}$. Define a new norm on $X$ by
\[\|x\|_1 = \sup \{x^*(x) : x^* \in B(F)\} \quad \mb{for } x
\in X\] 
Then $\|\cdot\|_1$ is a norm on $X$ with $F$ as a proper norming
subspace, and it follows from the proof of
\cite[Theorem~8.2]{GK} that this norm is equivalent to the
original norm. \qed

\Rem
In \cite{HL2} the authors showed that $X$ is reflexive if and
only if for any equivalent norm, $X$ is Hahn-Banach smooth and
has ANP-III. This was strengthened in \cite{br} to just
Hahn-Banach smooth (Corollary 2.5). The above is even stronger
result with even easier proof.
\end{rem}

Now we obtain some stability results for nicely smooth spaces.

\begin{thm} \label{l_p}
Let $\{X_\a\}_{\a\in \G}$ be a family of Banach spaces. Then $X
= \bigoplus_{\ell_p} X_\a$ $(1< p < \iy)$ is nicely smooth if
and only if for each $\a\in \G$, $X_\a$ is nicely smooth.
\end{thm}

\pf We will show that $C_X = \{0\}$ if and only if for every
$\a\in \G$, $C_{X_\a} = \{0\}$.

Now, $X = \bigoplus_{\ell_p} X_\a$ implies $X^{**} =
\bigoplus_{\ell_p} X^{**}_\a$, and $x^{**} \in C_X$ if and only
if
\beqa
\|x^{**} + \hat{x}\|_p & \geq & \|x\|_p \mb{ for all } x \in X\\
\Lglra \sum_{\a \in \G} \|x_\a^{**} + \hat{x}_\a\|^p & \geq &
\sum_{\a \in \G} \|x_\a\|^p \mb{ for all } x \in X
\eeqa

It is immediate that if for every $\a\in \G$, $x^{**}_\a \in
C_{X_\a}$, then $x^{**} \in C_X$. And hence, $C_X = \{0\}$
implies for every $\a\in \G$, $C_{X_\a} = \{0\}$.

Conversely, suppose for every $\a\in \G$, $C_{X_\a} = \{0\}$.
Let $x^{**} \in X^{**} \setminus \{0\}$. Let $\a_o \in \G$ be
such that $x^{**}_{\a_o} \neq 0$. Then $x^{**}_{\a_o} \notin
C_{X_{\a_o}}$. Hence, there exists $x_{\a_o} \in X_{\a_o}$ such
that $\|x^{**}_{\a_o} + \hat{x}_{\a_o}\| < \|x_{\a_o}\|$. Choose
$\e > 0$ such that $\|x^{**}_{\a_o} + \hat{x}_{\a_o}\|^p + \e <
\|x_{\a_o}\|^p$. Then there exists a finite $\G_o \ci \{\a \in
\G : x^{**}_\a \neq 0\}$ such that $\a_o \in \G_o$ and $\sum_{\a
\notin \G_o} \|x_\a^{**}\|^p < \e$. If $\a \in \G_o$, then
$x^{**}_\a \notin C_{X_\a}$. Hence, there exists $x_\a \in X_\a$
such that $\|x^{**}_\a + \hat{x}_\a\| < \|x_\a\|$. Define $y \in
X$ by
\[y_\a = \alto{x_\a} {\a \in \G_o}{0}\] 
Then we have,
\beqa
\|x^{**} + \hat{y}\|_p^p & = & \sum_{\a \in \G} \|x_\a^{**} +
\hat{y}_a\|^p \\ & = & \sum_{\a \in \G_o \above0pt \a \neq \a_0}
\|x_\a^{**} + \hat{x}_a\|^p + \|x_{\a_o}^{**} +
\hat{x}_{\a_o}\|^p + \sum_{\a \notin \G_o} \|x_\a^{**}\|^p \\
& < & \sum_{\a \in \G_o \above0pt \a \neq \a_0} \|x_a\|^p +
\|x_{\a_o}^{**} + \hat{x}_{\a_o}\|^p + \e\\
& < & \sum_{\a \in \G_o} \|x_\a\|^p = \|y\|_p^p
\eeqa
which shows that $x^{**} \notin C_X$. \qed

\Rem (a) The above argument also works for finite $\ell_1$ (or
$\ell_\iy$) sums and shows that if $X$ is the $\ell_1$ (or
$\ell_\iy$) sum of $\et{X}{n}$, then $X$ is nicely smooth if and
only if for every coordinate space $X_i$ is so.

However, if $\G$ is infinite, $X = \bigoplus_{\ell_1} X_\a$ is
never nicely smooth as $\bigoplus_{c_o} X^*_\a$ is a proper
norming subspace of $X^* = \bigoplus_{\ell_\iy} X^*_\a$.

A similar argument also shows that being nicely smooth is not stable
under infinite $\ell_\iy$ sums.

(b) Since Property $(II)$ is not preserved under finite $\ell_1$
sums \cite[Proposition 3.7]{br}, such sums produce examples of
nicely smooth spaces lacking Property $(II)$.
\end{rem}

We now show that being nicely smooth is stable under $c_o$ sums.

\begin{thm} \label{c_o}
Let $\{X_\a\}_{\a\in \G}$ be a family of Banach spaces. Then $X =
\bigoplus_{c_o} X_\a$ is nicely smooth if and only if for each
$\a\in \G$, $X_\a$ is nicely smooth.
\end{thm}

\pf As before, we will show that $C_X = \{0\}$ if and only if
for every $\a\in \G$, $C_{X_\a} = \{0\}$.

Necessity is similar to that in Theorem~\ref{l_p}.

Conversely, suppose for every $\a\in \G$, $C_{X_\a} = \{0\}$.
And let $x^{**} \in X^{**} \setminus \{0\}$. Let $\a_o \in \G$
be such that $x^{**}_{\a_o} \neq 0$. Then $x^{**}_{\a_o} \notin
C_{X_{\a_o}}$. Hence, there exists $x_{\a_o} \in X_{\a_o}$ such
that $\|x^{**}_{\a_o} + \hat{x}_{\a_o}\| < \|x_{\a_o}\|$.
Triangle inequality shows that for any $\la \geq 1$,
$\|x^{**}_{\a_o} + \la \hat{x}_{\a_o}\| < \|\la x_{\a_o}\|$.
Thus, replacing $x_{\a_o}$ by $\la x_{\a_o}$ for some $\la \geq
1$, if necessary, we may assume $\|x_{\a_o}\| >
\|x^{**}\|_\iy$.

Define $y \in X$ by
\[y_\a = \alto{x_{\a_o}} {\a = \a_o}{0}\] 
Then,
\[\|x^{**} + \hat{y}\|_\iy = \max\{\sup \{\|x^{**}_\a\|_{\a
\neq \a_o}\}, \|x^{**}_{\a_o} + \hat{x}_{\a_o}\|\} <
\|x_{\a_o}\| = \|y\|_\iy \] whence $x^{**} \notin C_X$. \qed

\begin{cor}
Being nicely smooth is not a three space property.
\end{cor}

\pf Let $X = c$, the space of all convergent sequences with the
sup norm. Recall that $c^* = \ell_1$ and that $\ell_1$ acts on
$c$ as
\[\langle \bm{a}, \bm{x} \rangle = a_0 \lim x_n + \sm{n=0}{\iy}
a_{n+1} x_n, \quad \bm{a} = \{a_n\}_{n=0}^{\iy} \in
\ell_1, \, \bm{x} = \{x_n\}_{n=0}^{\iy} \in c\]
It follows that $\{\bm{a} \in \ell_1 : a_0 = 0\}$ is a proper
norming subspace for $c$. 

Put $Y = c_o$. Then $Y$ is nicely smooth and dim$(X/Y) = 1$, so
that $X/Y$ is also nicely smooth. But, by above, $X$ is 
not nicely smooth. \qed

Recall that a closed subspace $M \ci X$ is said to be an
M-summand if there is a projection $P$ on $X$ with range $M$
such that $\|x\| = \max\{\|Px\|, \|x-Px\|\}$ for all $x \in X$.
An easy modification of the arguments of \cite[Proposition
2]{BR} shows that
\begin{prop}
If $Y$ is a $M$-summand in $X$ and $X$ has the BGP, then so does
$Y$.
\end{prop}

\begin{thm}
Let $X$ be a Banach space, $\m$ denote the Lebesgue measure on
$[0, 1]$ and $1<p< \iy$. \TFAE
\bla
\item $L_p(\m, X)$ has BGP.
\item $L_p(\m, X)$ is nicely smooth.
\item $X$ is nicely smooth and Asplund.
\el
\end{thm}

\pf Clearly $(a) \Ra (b)$.

$(b) \Ra (c)$. Since $L^q(\m, X^*)$ is always a norming subspace
of $L^p(\m, X)^*$, $\rec{p} + \rec{q} =1$, and they coincide if
and only if $X^*$ has the RNP with respect to $\m$ \cite[Chapter
IV]{DU}, $(b)$ implies $X^*$ has the RNP, or, $X$ is Asplund.
Also for any norming subspace $F \ci X^*$, $L^q(\m, F)$ is a
norming subspace of $L^p(\m, X)^*$. Hence, $(b)$ also implies
$X$ is nicely smooth.

$(c) \Ra (a)$. If $X$ is nicely smooth and Asplund, by
Corollary~\ref{ns}, $X^*$ is the closed linear span of the
w*-denting points of $B(X^*)$. And it suffices to show that
$L^p(\m, X)^* = L^q(\m, X^*)$ is the closed linear span of the
w*-denting points of $B(L^q(\m, X^*))$.

Let $F = \sm {i=1}{n} \a_i x^*_i \c_{A_i}$ with $x^*_i \in
S(X^*)$ for all $i = 1, 2, \ld, n$ be a simple function in
$S(L^q(\m, X^*))$. Let $\e > 0$. Now, for each $i = 1, 2, \ld,
n$, there exists $\la_{ik} \in \R$, and $x^*_{ik}$, w*-denting
points of $B(X^*)$, $k = 1, 2, \ld, N$, such that $\|x^*_i -
\sm{k=1}{N} \la_{ik} x^*_{ik}\| < \e$. For $k = 1, 2, \ld, N$.
Define
\[ F_k = \sm{i=1}{n} \a_i \la_{ik} x^*_{ik} \c_{A_i} \] 
Since each $x^*_{ik}$ is a w*-denting point of $B(X^*)$, for
each $k$, $F_k/\|F_k\|$ is a w*-denting point of $B(L^q(\m,
X^*))$ \cite[Lemma 10]{BR}. And,
\beqa 
\|F - \sm{k=1}{N} F_k\|_q^q &=& \|\sm {i=1}{n} \a_i x^*_i
\c_{A_i} - \sm{k=1}{N} \sm{i=1}{n} \a_i \la_{ik} x^*_{ik}
\c_{A_i}\|_q^q \\ &=&
\sm {i=1}{n} |\a_i|^q \|x^*_i - \sm{k=1}{N} \la_{ik}
x^*_{ik}\|^q \m(A_i) \\ & < &
\sm {i=1}{n} \e^q |\a_i|^q \m(A_i) \leq \e^q \|F\|^q_q \leq \e
\qquad \qed 
\eeqa 

The following results closely parallel the corresponding results
in \cite{br}. We include the proofs for completeness.

\begin{prop} \label{CKX}
Let $K$ be a compact Hausdorff space, then $C(K, X)$ is nicely
smooth if and only if $K$ is finite and $X$ is nicely smooth.
\end{prop}

\pf For a compact Hausdorff space $K$ and a Banach space $X$,
the set
\[A = \{\d(k) \ot x^* : k \in K, \, x^* \in S(X^*)\} \ci
B(C(K, X)^*)\]
is a norming set for $C(K, X)$. So, if $C(K, X)$ is nicely
smooth, $C(K, X)^* = \ov{span}(A)$. It follows that $K$ admits
no nonatomic measure, whence $K$ is scattered. Now, let $K'$
denote the set of isolated points of $K$. Then $K'$ is dense in
$K$, so, the set
\[A' = \{\d(k) \ot x^* : k \in K', \, x^* \in S(X^*)\}\]
is also norming. Thus, $C(K, X)^* = \ov{span}(A')$. But if $k
\in K \setminus K'$, then for any $x^* \in S(X^*)$, $\d(k)
\ot x^* \notin \ov{span}(A')$. Hence, $K = K'$, whence $K$
must be finite. And if $k_o \in K'$, $x \ra \chi_{\{k_o\}} x$ is
an isometric embedding of $X$ into $C(K, X)$ as an $M$-summand.
Thus, $X$ is nicely smooth.

The converse is immediate from Theorem~\ref{c_o}. \qed

\Rem (a) It is immediate that for $C(K)$ spaces Property $(II)$,
the BGP and being nicely smooth (indeed, any of the conditions
of Proposition~\ref{p1}) are equivalent, and are equivalent to
reflexivity. 

(b) It follows from the above and \cite[Theorem 3.9]{br} that
$C(K, X)$ has Property $(II)$ if and only if $K$ is finite and
$X$ has Property $(II)$.
\end{rem}

\begin{prop} \label{LX}
Let $X$ be a Banach space such that there exists a bounded net
$\{K_\a\}$ of compact operators such that $K_\a \lgra Id$ in the
weak operator topology. If ${\cal L}(X)$ is nicely smooth, then
$X$ is finite dimensional.
\end{prop}

\pf For $x\in X$, $x^*\in X^*$, let $x \ot x^*$ denote the
functional defined on ${\cal L}(X)$ by $(x \ot x^*)(T) =
x^*(T(x))$. Then $\|x\ot x^*\| = \|x\| \|x^*\|$. And, since
$\|T\| = \sup_{\|x^*\|= 1, \|x\|= 1} (x^*(T(x)) = \sup_{\|x^*\|=
1, \|x\|= 1} (x \ot x^*)(T)$, it follows that $A = \{x \ot
x^* : \|x^*\|= 1, \|x\|= 1\}$ is a norming set, and hence,
${\cal L}(X)^* = \ov{span}(A)$.

Claim : $K_\a \lgra Id$ weakly.

Since $\{K_{\a}\}$ is bounded, it now suffices to check that
$K_\a \lgra Id$ on $A$, i.e., to check $x^*(K_\a(x)) \lgra
x^*(x)$ for all $\|x\|=1$, $\|x^*\|= 1$.  But, $K_\a (x) \lgra
x$ weakly, hence the claim.

Thus, $Id$ is a compact operator, so that $X$ is finite
dimensional. \qed

\begin{prop}
If ${\cal L}(X^*)$ is nicely smooth, then $X$ is finite
dimensional.
\end{prop}

\pf Since ${\cal L}(X^*) = (X\bigotimes_\pi X^*)^*$ (i.e., the
projective tensor product of $X$ and $X^*$) is a nicely smooth
dual space, it is reflexive. But then $X$ and ${\cal L}(X)$ are
reflexive too. It follows from \cite[Theorem 2]{K} that $X$ must
be separable, and hence it admits a Marshukevich basis. The
corresponding finite rank projections produce a bounded sequence
of compact operators converging to $Id$ in the weak operator
topology. Now an argument similar to Proposition~\ref{LX} shows
that $X$ is finite dimensional.
\qed

\begin{prop} 
For a compact Hausdorff space $K$, ${\cal L}(X, C(K))$ is nicely
smooth if and only if ${\cal K}(X, C(K))$ is nicely smooth if
and only if $X$ is reflexive and $K$ is finite.
\end{prop}

\pf Suppose ${\cal L}(X, C(K))$ is nicely smooth. By definition
of the norm, $A = \{\d(k) \ot x : x \in B(X), k \in K\}$ is
a norming set for ${\cal L}(X, C(K))$, and hence, ${\cal L}(X,
C(K))^* = \ov{span}(A)$. It follows that ${\cal L}(X, C(K)) =
{\cal K}(X, C(K))$ and that ${\cal K}(X, C(K))$ is nicely
smooth.

Now, from the easily established identification, ${\cal K}(X,
C(K)) = C(K, X^*)$ and Proposition~\ref{CKX}, it follows that
${\cal K}(X, C(K))$ is nicely smooth if and only if $K$ is
finite and $X^*$ is nicely smooth, which, in turn, is equivalent
to $K$ is finite and $X$ is reflexive.

Also, if $K$ is finite, $C(K)$ is finite dimensional, so that
${\cal L}(X, C(K)) = {\cal K}(X, C(K))$. This completes the
proof. \qed

\bbibl
\small 
\setlength{\itemsep}{1pt}
\bibitem{BR} P.\ Bandyopadhyay and A.\ K.\ Roy, {\em Some
Stability Results for Banach Spaces with the Mazur Intersection
Property}, Indagatione Math.\ New Series {\bf 1}, No.\ 2,
(1990), 137--154.
\bibitem{br} Sudeshna Basu and TSSRK Rao, {\it Some Stability
Results for Asymptotic Norming Properties of Banach Spaces}, ISI
Tech.\ Report No.\ 1/96, (1996) (communicated).
\bibitem{CHL} D.\ Chen, Z.\ Hu and B.\ L.\ Lin, {\it Balls
intersection properties of Banach Spaces}, Bull.\ Austral.\
Math.\ Soc., {\bf 45}, (1992), 333--342.
\bibitem{CL} D.\ Chen and B.\ L.\ Lin, {\it Ball Separation
Properties in Banach spaces}, Preprint (1995).
\bibitem{CL1} D.\ Chen and B.\ L.\ Lin, {\it Ball Topology of
Banach spaces}, Preprint (1995).
\bibitem{DU} J.\ Diestel and J.\ J.\ Uhl, Jr., {\em Vector
Measures}, Math.\ Surveys {\bf 15}, Amer.\ Math.\ Soc.,
Providence, R.\ I.\ (1977).
\bibitem{G3} G.\ Godefroy, {\it Nicely Smooth Banach Spaces},
Longhorn Notes, The University of Texas at Austin, Functional
Analysis Seminar (1984--85), 117--124.
\bibitem{G4} G.\ Godefroy, {\it Existence and Uniqueness of
Isometric Preduals : A Survey}, Contemp.\ Math., {\bf 85}, (1989),
131--193.
\bibitem{GK} G.\ Godefroy and N.\ J.\ Kalton, {\it The ball
topology and its application}, Contemp.\ Math., {\bf 85}, (1989),
195--237. 
\bibitem{GS} G.\ Godefroy and P.\ D.\ Saphar, {\it Duality in
spaces of operators and Smooth norms on Banach Spaces}, Illinois
J.\ Math., {\bf 32}, (1988), 672--695.
\bibitem{HL2} Z.\ Hu and B.\ L.\ Lin, {\it Smoothness and
Asymptotic Norming Properties in Banach Spaces}, Bull.\
Austral.\ Math.\ Soc., {\bf 45}, (1992), 285--296.
\bibitem{K} N.\ J.\ Kalton, {\it Spaces of Compact Operators},
Math.\ Ann., {\bf 108}, (1974), 267--278.
\bibitem{R} H.\ Rosenthal, {\em On the Structure of Non-dentable
Closed Bounded Convex Sets}, Advances in Math.\ {\bf 70}, (1988),
1--58.
\ebibl
\ed

---------------------------end nice.tex----------------------------------------